\documentclass[12pt,a4paper,reqno]{amsart}
\usepackage{geometry}
\usepackage{graphicx}
\usepackage{placeins}
\usepackage[numbers]{natbib}

\newtheorem*{theorem*}{Theorem}

\theoremstyle{remark}
\newtheorem*{remark}{Remark}

\numberwithin{equation}{section}

\begin{document}
\title[Heuristics for discrete mean values]{A heuristic for discrete mean values of the derivatives of the Riemann zeta function}

\author{Christopher Hughes}
\address{Department of Mathematics, University of York, York, YO10 5DD, United Kingdom}
\email{christopher.hughes@york.ac.uk}
\author{Greg Martin}
\address{Department of Mathematics, University of British Columbia, Vancouver, British Columbia, V6T 1Z2, Canada}
\email{gerg@math.ubc.ca}
\author{Andrew Pearce-Crump}
\address{Department of Mathematics, University of York, York, YO10 5DD, United Kingdom}
\email{andrew.pearce-crump@york.ac.uk}

\begin{abstract}
Shanks conjectured that $\zeta ' (\rho)$, where $\rho$ ranges over non-trivial zeros of the Riemann zeta function, is real and positive in the mean. We present a history of this problem and its proof, including a generalization to all higher-order derivatives $\zeta^{(n)}(s)$, for which the sign of the mean alternatives between positive for odd~$n$ and negative for even~$n$. Furthermore, we give a simple heuristic that provides the leading term (including its sign) of the asymptotic formula for the average value of $\zeta^{(n)}(\rho)$.
\end{abstract}

\maketitle

\section{Introduction}\label{Sect:Intro}
Shanks' Conjecture (now a theorem), which dates from 1961, states that $\zeta '(\rho)$ is real and positive in the mean as $\rho = \beta + i \gamma$ ranges over non-trivial zeros of the Riemann zeta function $\zeta (s)$.
More recently, this assertion has been generalized to higher derivatives: On average, $ \zeta^{(n)}(\rho)$  is positive if $n$ is odd, but negative if $n$ is even.

The aim of this paper is to show how one can heuristically derive this Generalized Shanks' Conjecture, indeed in quantitative form, almost immediately from an explicit formula known as the Landau--Gonek Theorem.

The Riemann zeta function $\zeta(s)$ has infinitely many non-trivial zeros, which are numbers $\rho=\beta+i\gamma$ with $0<\beta<1$ that satisfy $\zeta(\rho)=0$. The locations of these non-trivial zeros, as we have known for a century and a half (since Riemann), are intimately connected with the distribution of prime numbers. While impressive calculations were done by hand over the years (culminating in Titchmarsh and Comrie's verification~\cite{Titchmarsh} that the first thousand zeros of $\zeta(s)$ lie on the critical line $\beta=\frac12$), the advent of electronic computers in the middle of the 20th century provided a significant boost to computations of the non-trivial zeros.

\begin{figure}[ht] \label{first derivative}
\centering
\includegraphics[height=6cm]{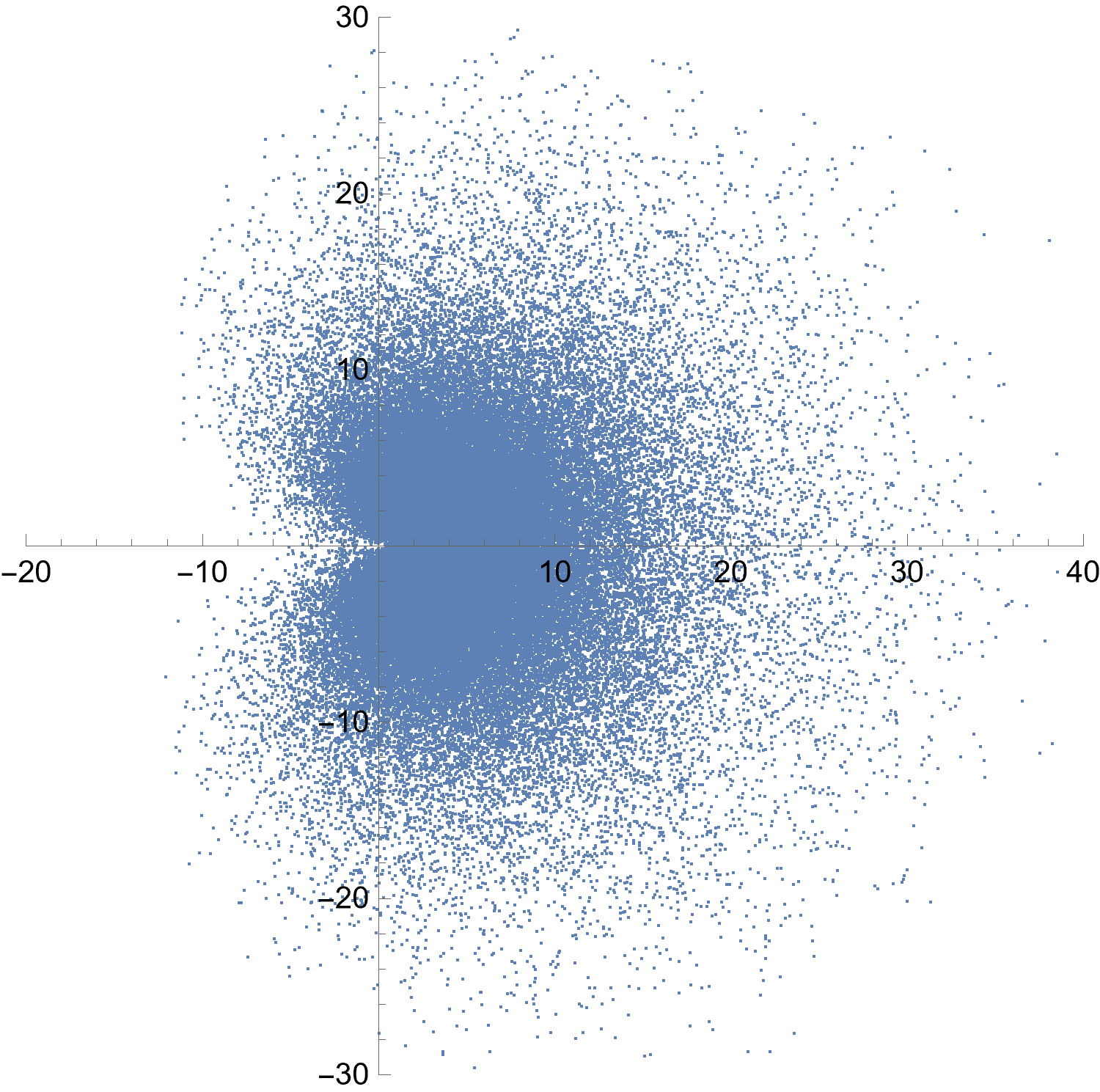}
\caption{Scatterplot of $\zeta'(\rho)$ for the first 100,000 zeros of the zeta function}
\end{figure}

Shanks first made his conjecture \cite{ShanksConj} when he was reviewing Haselgrove's tables \cite{Has60} of numerical values of the Riemann zeta function. He plotted the  graph of $t \mapsto \zeta ( \frac{1}{2} + it)$ and noticed that the way this curve approaches the origin,  mainly through the third and fourth quadrants, suggests the phase of $\zeta'(\frac{1}{2} + i \gamma)$ is close to zero in the mean and thus $\zeta'(\frac{1}{2} + i \gamma)$ is positive and real in the mean. This observation stands in contrast to the general behavior of the function $\zeta'(\frac12+it)$, whose mean value tends quickly to~$0$.
We present a different graph to the one Shanks created. Figure~\ref{first derivative} shows the graph of $\zeta ' (\rho)$ for the first 100,000 zeros of $\zeta (s)$, which clearly displays symmetry and a bias towards the positive side of the complex plane.

In 1985, Conrey, Ghosh and Gonek \cite{CGG88} were looking for a straightforward proof of the fact that there are infinitely many simple zeros of $\zeta (s)$. To this end, they used the Cauchy--Schwarz inequality to write
\begin{equation*}
    \bigg| \sum_{0 < \gamma \leq T} \zeta ' (\rho) \bigg| ^2 \leq  \left(\ \sideset{}{^*}\sum_{0 < \gamma \leq T} 1\right) \left(  \sum_{0 < \gamma \leq T} \left| \zeta ' (\rho) \right|^2 \right),
\end{equation*}
where the star on the middle sum denotes counting only the simple zeros. To complete the proof they used the asymptotic for the discrete second moment on the right-hand side, which was proved by Gonek \cite{Gonek} in 1984, who showed under the Riemann Hypothesis that
\[
\sum_{0 < \gamma \leq T} \left| \zeta ' (\rho) \right|^2 = \frac{T}{24\pi} (\log T)^4 + O\left(T (\log T)^3\right) ,
\]
and so they only needed to find the leading order behavior of the sum on the left-hand side. They finished their proof that there are infinitely many simple zeros of $\zeta (s)$ by showing that
\begin{equation}
\sum_{0 < \gamma \leq T} \zeta ' (\rho) = \frac{T}{4 \pi} (\log T)^2 + O (T \log T), \label{1stDeriv}
\end{equation}
which proved Shanks' conjecture as an additional benefit.

This approach shows the number of simple zeros is greater than some constant times $T$, but was not strong enough to show that a positive proportion of the non-trivial zeros of $\zeta(s)$ are simple (a fact that goes back to work of Levinson~\cite{Levinson}). However, later in 1998~\cite{CGG98}, they were able to prove, by modifying the method to include mollifiers, that at least $\frac{19}{27}$ of these non-trivial zeros are simple assuming the Riemann Hypothesis and the Generalized Lindel\"{o}f Hypothesis. In 2013, Bui and Heath-Brown~\cite{BHB13} removed the assumption of the Generalized Lindel\"{o}f Hypothesis through careful use of the generalized Vaughan identity.

Before we describe the generalization of Shanks' conjecture to higher derivatives, it is worth describing the ongoing research progress on the distribution of $\zeta'(\rho)$.
In 1994, Fujii~\cite{Fuj94} found explicit lower order terms for the asymptotic formula, Equation~\eqref{1stDeriv} (in \cite{FujiiDist}, he later corrected a slight error in the lowest-order coefficient). The full correct asymptotic is given~by
\begin{equation*}
    \sum_{0 < \gamma \leq T} \zeta ' (\rho) = \frac{T}{4 \pi} \biggl(\log \frac{T}{2\pi}\biggr)^2 + (-1+C_0) \frac{T}{2 \pi} \log \frac{T}{2\pi} + (1-C_0-C_0^2 + 3C_1) \frac{T}{2\pi} + E(T),
\end{equation*}
where $C_0$ and $C_1$ are coefficients in the Laurent expansion of $\zeta (s)$ about $s=1$. Unconditional effective bounds for $E(T)$ are known, and a savings of a power of~$T$ can be shown under the Riemann Hypothesis. The best known unconditional error term can be found in~\cite{HugPC22}. In 2021, Kobayashi \cite{Kob21a} gave a similar result for Dirichlet $L$-functions.

In 2010, Trudgian  \cite{Trudgian2010} gave an alternative proof of Shanks' conjecture, based on Shanks' observation of the connection between $\arg \zeta'(\rho)$ and the so-called Gram's Law. Although care needs to be taken in the definition of the argument (a priori it is only defined up to a multiple of $2\pi$), Trudgian was able to show that
\[
\sum_{0 < \gamma \leq T} \arg \zeta'(\rho) \ll_\varepsilon T^\varepsilon
\]
for every $\varepsilon>0$.

We are now ready to explore the generalization of Shanks' Conjecture to $\zeta^{(n)} (s)$, the $n$th derivative of $\zeta (s)$, for every positive integer~$n$. This extension, which we call the Generalized Shanks' Conjecture, states that $\zeta^{(n)}(\rho)$ is real and positive in the mean if $n$ is even, while $\zeta^{(n)}(\rho)$ is real and negative in the mean if $n$ is odd.

\begin{figure}[ht] \label{higher derivatives}
\begin{center}
     \begin{tabular}{p{0.45\textwidth}  p{0.45\textwidth}}
    \includegraphics[width=0.4\textwidth]{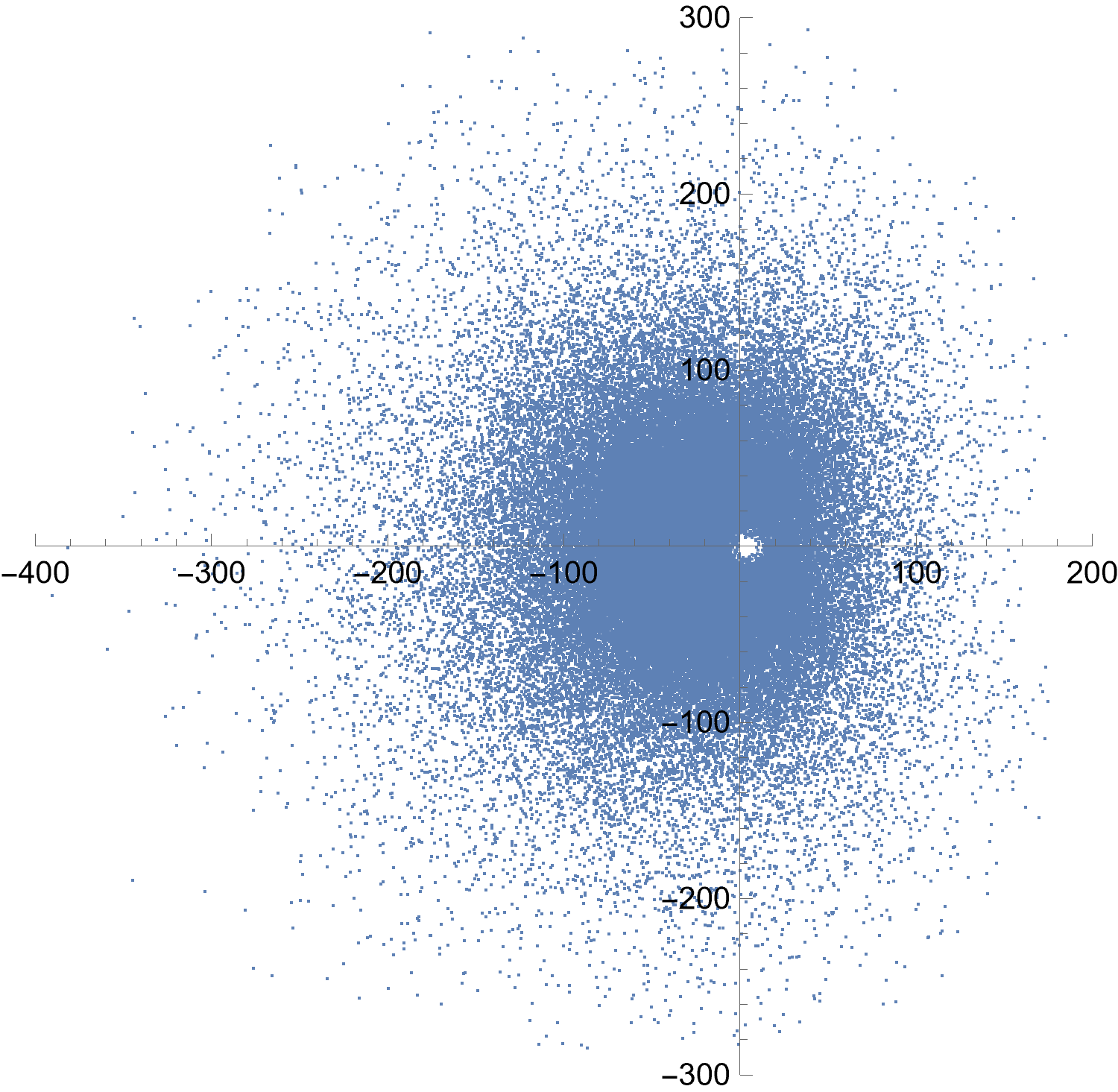}
      &
      \includegraphics[width=0.4\textwidth]{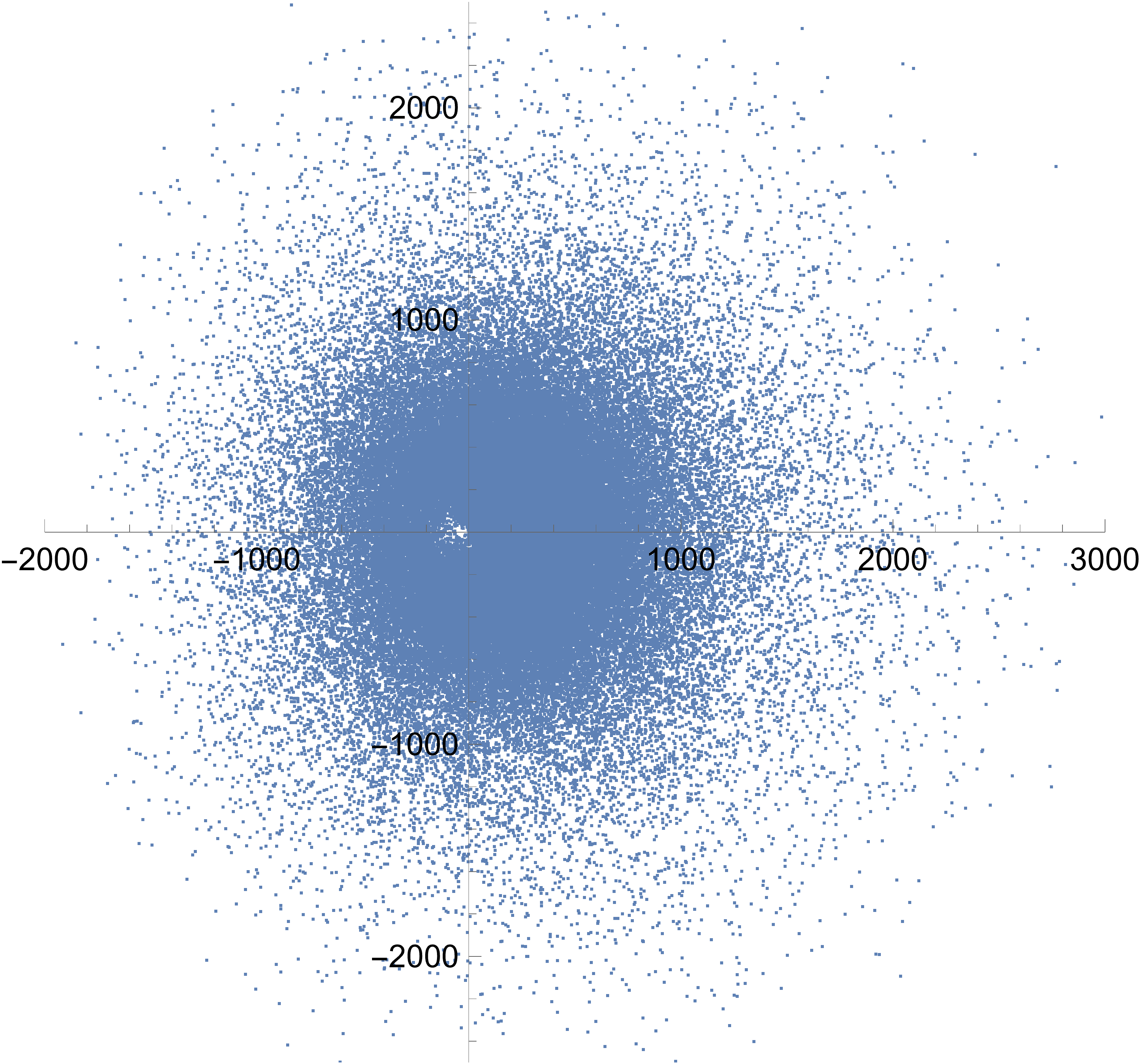}

      \end{tabular}
      \end{center}
\caption{Scatterplot of $\zeta^{(n)}(\rho)$ for $n=2$ and $n=3$ for the first 100,000 zeros of the zeta function}
\end{figure}

In 2011, using ideas similar to those in Conrey, Ghosh and Gonek's proof of Shanks' original conjecture, Kaptan, Karabulut and Y{\i}ld{\i}r{\i}m \cite{KKY11} found the main term for the summatory function of $\zeta^{(n)}(\rho)$. Specifically, they showed that
\begin{equation}\label{nthderiv}
\sum_{0 < \gamma \leq T} \zeta^{(n)} (\rho) =
  \frac{(-1)^{n+1}}{n+1} \frac{T}{2 \pi} (\log T)^{n+1} + O \bigl(T (\log T)^n \bigr),
\end{equation}
which clearly implies the Generalized Shanks' Conjecture as stated above. They have also proven a similar result for derivatives of Dirichlet $L$-functions, with the same leading-order asymptotic term. Note that Shanks' Conjecture was first conjectured and later proved, in the usual order of events, but that we still call the established assertion Shanks' ``Conjecture''. This nomenclature leads to an unusual situation where the extension is reasonably labeled ``Shanks' Generalized Conjecture'' even though its first appearance in the literature was as a proven result!

In 2022, the first and third author~\cite{HugPC22} investigated the Generalized Shanks' Conjecture and found all the lower order terms for the asymptotic formula. The result can be stated in the form
\[
\sum_{0 < \gamma \leq T} \zeta^{(n)} (\rho) =
  \frac{(-1)^{n+1}}{n+1} \frac{T}{2 \pi} \biggl( \log \frac{T}{2 \pi} \biggr)^{n+1}  +\frac{T}{2\pi} \mathcal{P}_n \biggl(\log \frac{T}{2 \pi} \biggr) + E_n(T),
\]
where $\mathcal{P}_n (x)$ is a polynomial of degree $n$ described explicitly in \cite{HugPC22}; as with Fujii's result, unconditional effective bounds for the error term $E_n(T)$ are known, and a power savings in~$T$ can be shown under the Riemann Hypothesis. Furthermore, the third author~\cite{Pearce-Crump21} has also extended Fujii's combination (mentioned earlier) of $\zeta ' (s)$ and the Landau--Gonek Theorem to all derivatives of $\zeta (s)$.

All of the results we have described have complicated proofs, and it is unclear whether any of them provide an intuitive explanation for the alternating signs in the Generalized Shanks' Conjecture, or indeed of the positive sign in Shanks' original conjecture. We now present a simple heuristic that explains these signs, and indeed recovers (non-rigorously) the leading-order term in the asymptotic formula given in Equation~\eqref{nthderiv}.

\section{Heuristic for the Riemann Zeta Function}

We begin with some preliminaries concerning the Riemann zeta function, which for $\Re (s) > 1$ can be written as the convergent Dirichlet series
\begin{equation*}
    \zeta (s) = \sum_{m=1}^{\infty} \frac{1}{m^s}.
\end{equation*}
Since Dirichlet series converge locally uniformly we may differentiate term by term; and since the derivative of $\frac1{m^s} = e^{-s\log m}$ is $({-}\log m)m^{-s}$, we see that
\begin{equation*}
    \zeta^{(n)} (s) = (-1)^n \sum_{m=1}^{\infty} \frac{(\log m)^n}{m^s}
\end{equation*}
for $\Re (s) > 1$. This formula can be extended beyond the range of absolute convergence in myriad ways. One such extension (see the appendix for a proof) states that uniformly for $\sigma\geq\sigma_0>0$ and $t>2$,
\begin{equation}
\zeta^{(n)} (s) = (-1)^n \sum_{m \leq t} \frac{(\log m)^n}{m^s} + O( t^{-\sigma}(\log t)^n) .\label{approxform}
\end{equation}

We state a corollary of the Landau--Gonek Theorem, found in \cite{Gon93}, concerning sums of the form $X^\rho$ over the non-trivial zeros of the Riemann zeta function. This corollary is the main result that the heuristic relies upon. We call this result the Landau--Gonek Theorem for the sake of brevity in this paper, but note that it is given in a more general form in the cited papers. It was proved by Landau \cite{Landau11} for fixed $X$ and made uniform by Gonek \cite{Gon85,Gon93}. We note that while our quoted result requires the assumption of the Riemann Hypothesis, the general result does not.

\begin{theorem*}[Landau--Gonek Theorem]
Under the Riemann Hypothesis, for $T>1, m \in \mathbb{N}$ with $m\geq 2$ and $\rho$ a non-trivial zero of the Riemann zeta function $\zeta (s)$,
\begin{equation*}
\sum_{0< \gamma \leq T} m^{-\rho} = - \frac{T}{2 \pi} \frac{ \Lambda (m)}{m} + O(\log (2mT) \log\log (3m)),
\end{equation*}
where  $\Lambda (m)$ is the von Mangoldt function, which therefore only contributes when $m$ is a prime power.
\end{theorem*}

\subsection{The Heuristic}

We now state the short heuristic argument which yields the Generalized Shanks' Conjecture as stated in Section~\ref{Sect:Intro}. We will mainly ignore error terms so the equal signs are not `true' equalities but show the general argument. (But see Section~\ref{sect:ErrorTerms}.)

\begin{remark}
As we were writing this paper, we became aware that the heuristic presented in this paper can be found in \cite{Gon93} in a different context, where it can be applied rigorously.
\end{remark}

Throughout this section we write $\rho = \beta + i \gamma$ for a non-trivial zero of $\zeta (s)$. Using the approximate formula for $\zeta (s)$ given in Equation~\eqref{approxform} and the Landau--Gonek Theorem, we may write for $n \geq 1$,
\begin{align}
    \sum_{0 < \gamma \leq T} \zeta^{(n)} (\rho) & \approx (-1)^n \sum_{0 < \gamma \leq T} \sum_{m \leq \gamma} (\log m)^n m^{-\rho} \label{eq:ApproxSumZeros}\\
    & = (-1)^n \sum_{m\leq T} (\log m)^n \sum_{m < \gamma \leq T} m^{-\rho} \notag\\
    & \approx (-1)^{n+1} \frac{T}{2 \pi}  \sum_{m\leq T} \frac{ (\log m)^n \Lambda (m)}{m} - (-1)^{n+1} \frac{1}{2\pi} \sum_{m\leq T} (\log m)^n \Lambda(m). \notag
\end{align}
To finish the heuristic, we need to sum the series in the last line. By Chebyshev's Theorem~\cite{Apost},
\begin{equation*}
    C(x) = \sum_{m \leq x} \frac{\Lambda (m)}{m} = \log x + O(1),
\end{equation*}
so by partial summation, we have
\begin{align*}
\sum_{m \leq T} \frac{ (\log m)^n \Lambda (m)}{m} & = C(T) (\log T)^n - n  \int_{1}^{T} \frac {C(x) (\log x)^{n-1}}{x} \ dx \\
& = \frac{1}{n+1} (\log T)^{n+1} + O\bigl((\log T)^{n}\bigr).
\end{align*}
and similarly
\[
\sum_{m\leq T} (\log m)^n \Lambda(m) = O\bigl(T (\log T)^n \bigr).
\]
Combining this with our argument above, we have
\begin{align*}
    \sum_{0 < \gamma \leq T} \zeta^{(n)} (\rho) & = (-1)^{n+1} \frac{T}{2 \pi} \frac{1}{n+1} (\log T)^{n+1} + O(T (\log T)^{n}),
\end{align*}
as $T \rightarrow \infty$, which is the leading order the asymptotic result for the Generalized Shanks' Conjecture.

\begin{remark}
Note that this result assumes $n\geq 1$. When $n=0$ the LHS is trivially zero, whilst the RHS, as written, is not. This is explained by noticing that for $n\geq 1$ the $m=1$ term in the sum in Equation~\eqref{eq:ApproxSumZeros} is not present, since $\log(1)^n = 0$. However, if $n=0$ this term is present and is not accounted for in the Landau--Gonek Theorem, which holds for $m\geq 2$ only. The $m=1$ term in that case clearly contributes $N(T)$, which perfectly cancels the $-\frac{T}{2\pi}\log T + O(T)$ term coming from the calculation given above.
\end{remark}

\begin{remark}
This argument can be applied \emph{mutatis mutandis} for Dirichlet $L$-functions, yielding exactly the same leading-order behavior as for the Riemann zeta function. This result was already known \cite{KKY11,Kob21a}.
\end{remark}

\subsection{The Error Terms} \label{sect:ErrorTerms}

In the previous section, in the last part of Equation~\eqref{eq:ApproxSumZeros} we ignored the error terms coming from the Landau--Gonek Theorem. Including them one can see that they contribute
\[
\ll \sum_{m\leq T} (\log m)^n \log\log(3m) \log T \ll T (\log T)^{n+1} \log\log T,
\]
which dominates the main term! However, these are worst-case point-wise estimates and take no account of any potential cancellation when averaged, so it is likely that when summed over $m$ the true error is smaller and the main term is correct. (Indeed that is what is proved, by different methods, in \cite{CGG88, KKY11} for $n=1$ and $n\geq 2$ respectively).

The point is that whilst this method must remain a heuristic, it is a much quicker approach to find the main term in the overage of $\zeta^{(n)}(\rho)$ and gives some sense of why the Generalized Shanks' Conjecture holds true, that is why the mean of $\zeta^{(n)}(\rho)$  is real and positive / negative.

\begin{remark}
    In \cite{HugPC22}, the first and third authors prove that the true value of the mean of $\zeta^{(n)}(\rho)$ (that is, the main term and all the lower-order terms up to a power saving) is given by
    \[
\sum_{\substack{\ell , m \\ \ell m \leq \frac{T}{2\pi}}} \Lambda(m) (\log m)^n ,
    \]
    which shows that the heuristic in the previous section can only yield the main term since the last line of Equation~\eqref{eq:ApproxSumZeros} differs from the true value by $O(T(\log T)^n)$.
\end{remark}

\vskip20pt\noindent {\bf Acknowledgements.}

The authors thank the anonymous referee for several helpful comments.
The second author suggested the heuristic in this paper at the PIMS CRG `Moments of $L$-functions' workshop in Prince George, Canada to the third author. The third author would like to thank the group, the LMS and the University of York for their generous funding.

\appendix

\section{Appendix: Proof of Equation~\eqref{approxform}}

We sketch the proof for Equation~\eqref{approxform}, using an adaptation of the method found in Titchmarsh~\cite[pages~74--77]{Titchmarsh86}. Equation (4.11.2) from that book shows that for $\sigma>0$ and $N \in \mathbb{N}$,
\[
\zeta(s) = \sum_{m=1}^N \frac{1}{m^s} - \frac{N^{1-s}}{1-s} + s\int_N^\infty \frac{[u]-u+1/2}{u^{s+1}} \ du - \frac12 N^{-s}.
\]
If we differentiate this equality $n$ times with respect to $s$, we get
\begin{multline}\label{eq}
\zeta^{(n)}(s) = (-1)^n \sum_{m=1}^N \frac{(\log m)^n}{m^s} - \frac{d^n}{ds^n} \left[ \frac{N^{1-s}}{1-s} \right] - \frac12 (-\log N)^n N^{-s} \\
+  n \int_N^\infty \frac{([u]-u+1/2)(-\log u)^{n-1}}{u^{s+1}} \ du +  s \int_N^\infty \frac{([u]-u+1/2) (-\log u)^{n}}{u^{s+1}} \ du.
\end{multline}
When $\sigma>\sigma_0$ for a fixed $\sigma_0>0$, the last three terms can be uniformly estimated as
\begin{align*}
- \frac12 (-\log N)^n N^{-s} &\ll \frac{(\log N)^n}{N^{\sigma_0}} \\
\int_N^\infty \frac{([u]-u+1/2)(-\log u)^{n-1}}{u^{s+1}} du  &\ll \frac{(\log N)^{n-1}}{N^{\sigma_0}} \\
s \int_N^\infty \frac{([u]-u+1/2)(-\log u)^{n}}{u^{s+1}} du  &\ll |s| \frac{(\log N)^n} {N^{\sigma_0}}.
\end{align*}
Moreover, if we define $g(m) = \frac{(\log m)^n}{m^\sigma}$ and $f(m) = -\frac{ t \log m }{2\pi}$,
we may apply~\cite[Lemma~4.10]{Titchmarsh86} to get
\begin{align*}
(-1)^n \sum_{t < m \leq N} \frac{(\log m)^n}{m^s} &= (-1)^n \sum_{t < m \leq N} g(m) e^{2\pi i f(m)}\\
&= (-1)^n \int_t^N g(u) e^{2\pi i f(u)} du + O\left( \frac{(\log t)^n}{t^\sigma} \right) \\
    &=  \int_t^N \frac{(-\log u)^n}{u^s} du + O\left( \frac{(\log t)^n}{t^\sigma} \right).
\end{align*}

On the other hand,
\[
 \int_t^N \frac{(-\log u)^n}{u^s} du = \frac{d^n}{ds^n} \left[ \frac{N^{1-s}}{1-s} \right] - \frac{d^n}{ds^n} \left[ \frac{t^{1-s}}{1-s} \right]
\]
(to see this, note that for $n=0$ the integral is directly calculable, and for $n>1$ simply differentiate both sides with respect to $s$ the appropriate number of times).
If we restrict ourselves to $t>2$, we obtain
\[
\frac{d^n}{ds^n} \left[ \frac{t^{1-s}}{1-s} \right] \ll \frac{t^{1-\sigma}(\log t)^n}{|1-\sigma + i t|} \ll t^{-\sigma}(\log t)^n .
\]
Therefore, plugging all this back into Equation~\eqref{eq}, we have shown that when $\sigma>\sigma_0$ for a fixed $\sigma_0>0$,
\[
\zeta^{(n)}(s) = (-1)^n \sum_{m \leq t} \frac{(\log m)^n}{m^s} + O\left( |s| \frac{(\log N)^n} {N^{\sigma_0}} \right) + O\left( \frac{(\log t)^n}{t^\sigma} \right);
\]
letting $N\to\infty$ completes the derivation of Equation~\eqref{approxform}.

\end{document}